\documentclass[11pt]{article}

\usepackage{amsmath,amssymb,latexsym,amsfonts,mathrsfs,graphicx,verbatim}

\usepackage{graphicx}
\usepackage{citesort}
\usepackage{setspace}
\usepackage{program}
\usepackage{subfig}

\renewcommand{\quad}{$~~~\;\;\;$}
\newtheorem{theorem}{Theorem}

\newtheorem{proposition}{Proposition}

\newtheorem{example}{Example}

\newcommand{\transp}{{^{\rm T}}}

\renewcommand{\R}{\mathbb{R}}
\newcommand{\matr}[1]{\begin{bmatrix} #1 \end{bmatrix}}    
\def\transp{^{\rm T}}

\newcommand{\Cond}{\mathcal C}
\newcommand{\Gr}{\mbox{\sf Gr}}

\providecommand{\newoperator}[3]{%
  \newcommand*{#1}{\mathop{#2}#3}}

\newoperator{\argmax}{\mathsf{argmax}}{}
\newoperator{\argmin}{\mathrm{argmin}}{}

\newcommand{\dmin}{\displaystyle\min}
\newcommand{\dmax}{\displaystyle\max}

\renewcommand{\span}{\mbox{\rm span}}

\author{ Javier Pe\~na\thanks{Tepper School of Business, Carnegie
Mellon University, USA, {\tt jfp@andrew.cmu.edu}}
\and
Vera Roshchina\thanks{Collaborative Research Network, University of Ballarat, AUSTRALIA, {\tt vroshchina@ballarat.edu.au}}
\and Negar Soheili\thanks{Tepper School of Business,
Carnegie Mellon University, USA, {\tt nsoheili@andrew.cmu.edu}}
}
\title{Some preconditioners for systems of linear inequalities}

\begin{document}
\maketitle
\abstract{We show that a combination of two simple preprocessing steps would generally improve the conditioning of a homogeneous system of linear inequalities.  Our approach is based on a comparison among three different but related notions of conditioning for linear inequalities.}  

\section{Introduction}
Condition numbers play an important role in numerical analysis.  The condition number of a problem is a key parameter in the computational complexity of iterative algorithms as well as in issues of numerical stability.  A related challenge of paramount importance is to {\em precondition} a given problem instance, that is, perform some kind of data preprocessing to transform a given problem instance into an equivalent one with that is better conditioned.  Preconditioning has been extensively studied in numerical linear algebra~\cite{GoluV96,TrefB97} and is an integral part of the computational implementation of numerous algorithms for solving linear systems of equations.  In the more general optimization context, the task of designing preconditioners has been studied by Epelman and Freund~\cite{EpelF02} as well as by Belloni and Freund~\cite{BellF09}.  In a similar spirit, we propose two simple preconditioning procedures for homogeneous systems of linear inequalities.  As we formalize in the sequel, these procedures lead to improvements in three types of condition numbers for homogeneous systems of linear inequalities, namely Renegar's~\cite{Rene95a}, Goffin-Cucker-Cheung's~\cite{cheung2001new,goffin1980relaxation}, and the Grassmann condition number~\cite{BellF09b,amelunxen2012coordinate}.  Both Renegar's and Goffin-Cheung-Chucker's condition numbers are key parameters in the analyses of algorithmic schemes and other numerical properties of constraint systems~\cite{EpelF00,FreuV99b,Rene95a,Rene95b,PenaR00,SoheP12}.  The more recently developed Grassmann condition number is  especially well-suited for probabilistic analysis~\cite{AmelB13}. 

We recall the definition of the above condition numbers in Section~\ref{sec.cnumbers} below.  These condition numbers quantify certain properties associated to the homogenous systems of inequalities \eqref{primal} and \eqref{dual} defined by a matrix $A\in \R^{m\times n},$ where $m \le n$.  Renegar's condition number is defined in terms of the distance from $A$ to a set of {\em ill-posed} instances in the space $\R^{m\times n}$.  Goffin-Cucker-Cheung's condition number is defined in terms of the {\em best conditioned} solution to the system defined by $A$. Alternatively, it can be seen as the reciprocal of a measure of thickness of the cone of feasible solutions to the system defined by $A$.  Goffin-Cucker-Cheung's condition number is intrinsic to certain geometry in the column space of $A$.  The more recent Grassmann condition number, introduced in a special case by Belloni and Freund~\cite{BellF09b} and subsequently generalized by Amelunxen-Burgisser~\cite{amelunxen2012coordinate} is based on the projection distance from the linear subspace spanned by the rows of $A$ to a set of {\em ill-posed} subspaces in $\R^n$.  The Grassmann condition number is intrinsic to certain geometry in the row space of $A$.

The Goffin-Cucker-Cheung's condition number and the Grassmann condition numbers are respectively invariant under column scaling and elementary row operations on the matrix $A$ respectively.  For suitable choices of norms, each of them is also always smaller than Renegar's condition number (see \eqref{eqn: CRvsGCC} and \eqref{eqn: CRvsAB} below).  We observe that these properties have an interesting parallel and naturally suggest two preconditioning procedures, namely {\em column normalization} and {\em row balancing.}  Our main results, presented in Section~\ref{sec.precond}, discuss some interesting properties of these two preconditioning procedures.  In particular, we show that a combination of them would improve the values of the three condition numbers.

\section{Condition numbers and their relations}\label{sec.cnumbers}
Given a matrix $A\in\R^{m\times n}$ with $m \le n$, consider the homogeneous feasibility problem\begin{equation}\label{primal}
Ax = 0,\hspace{5pt} x \geq 0, \hspace{5pt} x\neq 0,
\end{equation}
and its alternative 
\begin{equation}\label{dual}
A\transp y \ge 0, \;  y \neq 0.
\end{equation}
Let $F_P$ and $F_D$ denote the sets of matrices $A$ such that \eqref{primal} and \eqref{dual} are respectively feasible. The sets $F_P$ and $F_D$ are closed and $F_P \cup F_D = \R^{m\times n}$. 
The set $\Sigma := F_P\cap F_D$ is the set of {\em ill-posed} matrices.  For a given $A\in \Sigma$,   arbitrary small perturbations on $A$ can lead to a change with respect to the feasibility of \eqref{primal} and \eqref{dual}.

\bigskip

Renegar's condition number is defined as 
\[
\Cond_R(A) := \dfrac{\|A\|}{\min\{\|A-A'\|: A'\in \Sigma\}}.
\]
Here $\|A\|$ denotes the operator norm of $A$ induced by a given choice of norms in $\R^n$ and $\R^m$. 
When the Euclidean norms are used in both $\R^n$ and $\R^m$, we shall write  $\Cond^{2,2}_R(A)$ for $\Cond_R(A)$.  On the other hand, when the one-norm is used in $\R^n$ and the Euclidean norm is used in $\R^m$ we shall write $\Cond^{1,2}_R(A)$ for $\Cond_R(A)$.  These are the two cases we will consider in the sequel.

\bigskip

Assume $A =  \matr{a_1 & \cdots & a_n}\in\R^{m\times n}$ with $a_i\ne 0$ for $i=1,\dots,n$. Goffin-Cheung-Cucker's condition number is defined as
\[
\Cond_{GCC}(A) := \dmin_{\|y\|_2=1} \dmax_{i\in\left\{1,\ldots,n\right\}} \dfrac{\|a_i\|_2}{a_i\transp y}.
\] 
Here $\|\cdot\|_2$ denotes the Euclidean norm in $\R^m$.  
The quantity $1/\Cond_{GCC}(A)$ is the Euclidean distance from the origin to the boundary of the convex hull of the set of normalized vectors $\left\{\frac{a_i}{\|a_i\|_2}, i=1,\dots,n\right\}.$  Furthermore, when $A \in F_D$, this distance coincides with the thickness of the cone of feasible solutions to \eqref{dual}.

\bigskip

Let $\Gr_{n,m}$ denote the set of $m$-dimensional linear subspaces in $\R^n$, that is, the $m$-dimensional Grassmann manifold in $\R^n$.  Let
$$P_m:=\{W\in \Gr_{n,m}: W^\perp \cap \R^n_+ \ne \{0\}\},$$ 
$$D_m:=\{W\in \Gr_{n,m}: W \cap \R^n_+ \ne \{0\}\},$$ 
and $$\Sigma_m=P_m\cap D_m = \left\{W\in \Gr_{m,n}: W\cap \R^n_+\neq \left\{0\right\}, W\cap \mbox{int}( \R^n_+) = \emptyset \right\}.$$  

The sets  $P_m, D_m,$ and $\Sigma_m$ are analogous to the sets $F_P, F_D,$ and $\Sigma$ respectively:  If $A\in\R^{m\times n}$ is full-rank then $A\in F_P \Leftrightarrow \span(A\transp) \in P_m$ and   $A\in F_D \Leftrightarrow \span(A\transp) \in D_m$.  In particular, $A\in \Sigma \Leftrightarrow \span(A\transp) \in \Sigma_m.$  The Grassmann condition number is defined as~\cite{amelunxen2012coordinate}:
\[
\Cond_{Gr}(A) := \dfrac{1}{\min\{d(\span(A\transp),W):W\in\Sigma_m\}}.
\]
Here $d(W_1,W_2) = \sigma_{\max}(\Pi_{W_1} - \Pi_{W_2})$, where $\Pi_{W_i}$ denotes the orthogonal projection onto $W_i$ for $i=1,2.$ In the above expression and in the sequel, $\sigma_{\max}(M)$ and 
$\sigma_{\min}(M)$ denote the largest and smallest singular values of the matrix $M$ respectively.

\bigskip

We next recall some key properties of the condition numbers $\Cond^{1,2}_R,$ $\Cond^{2,2}_R,$ $\Cond_{GCC}$ and $\Cond_{Gr}$.  Throughout the remaining part of the paper assume $A = \matr{a_1 & \cdots & a_n}\in\R^{m\times n}$ is full-rank and $a_i\ne 0$ for $i=1,\dots,n$.  From the relationship between the one-norm and the Euclidean norm, it readily follows that
\begin{equation}\label{ineq.norms}
\frac{\Cond^{2,2}_R(A)}{\sqrt{n}} \le \Cond^{1,2}_R(A) \le \sqrt{n}\Cond^{2,2}_R(A),
\end{equation}
The following property is a consequence of \cite[Theorem 1]{cheung2001new}:
\begin{equation}\label{eqn: CRvsGCC}
\Cond_{GCC}(A) \leq \Cond^{1,2}_R(A) \leq \dfrac{\dmax_{i=1,\dots,n} \|a_i\|_2}{\dmin_{i=1,\dots,n}\|a_i\|_2}\Cond_{GCC}(A).
\end{equation}
The following property was established in \cite[Theorem 1.4]{amelunxen2012coordinate}:
\begin{equation}\label{eqn: CRvsAB}
\Cond_{Gr}(A) \leq \Cond^{2,2}_R(A) \leq \dfrac{\sigma_{\max}(A)}{\sigma_{\min}(A)}\Cond_{Gr}(A).
\end{equation}
The inequalities \eqref{eqn: CRvsGCC} and \eqref{eqn: CRvsAB} reveal an interesting parallel between the pairs  $\Cond_{GCC}(A), \Cond^{1,2}_R(A)$ and $\Cond_{Gr}(A), \Cond^{2,2}_R(A)$.  This parallel becomes especially striking by observing that 
the fractions  in the right hand sides of \eqref{eqn: CRvsGCC} and \eqref{eqn: CRvsAB} can be written respectively as 
\[
\frac{\dmax_{i=1,\dots,n}\|a_i\|_2}{\dmin_{i=1,\dots,n}\|a_i\|_2}  = \frac{\dmax\{\|Ax\|_2: \|x\|_1 = 1\}}{\min\{\|Ax\|_2: \|x\|_1 = 1\}},
\]
and
\[
\frac{\sigma_{\max}(A)}{\sigma_{\min}(A)} = \frac{\max\{\|A\transp y\|_2: \|y\|_2 = 1\}}{
\min\{\|A\transp y\|_2: \|y\|_2 = 1\}}.
\]

\section{Preconditioning via normalization and balancing}\label{sec.precond}

Inequalities \eqref{eqn: CRvsGCC} and \eqref{eqn: CRvsAB} suggest two preconditioning procedures for improving Renegar's condition number $\Cond_R(A)$.  The first one is to {\em normalize,} that is, scale the columns of $A$ so that the preconditioned matrix $\tilde A$ satisfies $\dmax_{i=1,\dots,n} \|\tilde a_i\| = \dmin_{i=1,\dots,n} \|\tilde a_i\|$.  The second one is to {\em balance,} that is, apply an orthogonalization procedure to the rows of $A$ such as Gram-Schmidt or QR-factorization so that the preconditioned matrix $\tilde A$ satisfies $\sigma_{\max}(\tilde A)=\sigma_{\min}(\tilde A).$  
Observe that if the initial matrix $A$ is full rank and has non-zero columns, these properties are preserved by each of the above two preconditioning procedures.
The following proposition formally states some properties of these two procedures.  
\begin{proposition}\label{prop.precond} Assume $A = \matr{a_1 & \cdots & a_n}\in\R^{m\times n}$ is full-rank and $a_i\ne 0$ for $i=1,\dots,n$. 
\begin{description}
\item[(a)] Let $\tilde A$ be the matrix obtained after normalizing the columns of $A$, that is, $\tilde a_i = \frac{a_i}{\|a_i\|}, \; i=1,\dots,n.$  Then
\[
\Cond^{1,2}_R(\tilde A) =  \Cond_{GCC}(\tilde A) = \Cond_{GCC}(A).
\]
This transformation is optimal in the following sense
\[
\Cond^{1,2}_R(\tilde A) = \min\{\Cond^{1,2}_R(AD): D \mbox{ is diagonal positive definite}\}.
\]
\item[(b)] Let $\tilde A$ be the matrix obtained after balancing the rows of $A$, that is, $\tilde A = Q\transp$, where $QR = A\transp$ with $Q\in\R^{n\times m},R\in\R^{m\times m}$ is the QR-decomposition of $A\transp$.
Then
\[
\Cond^{2,2}_R(\tilde A)= \Cond_{Gr}(\tilde A) =\Cond_{Gr}(A).
\]
This transformation is optimal in the following sense
\[
\Cond^{2,2}_R(\tilde A) = \min\{\Cond^{2,2}_R(PA): P\in \R^{m\times m} \mbox{ is non-singular}\}.
\]
\end{description}
\end{proposition}
\proof
Parts (a) and (b) follow respectively from \eqref{eqn: CRvsGCC} and \eqref{eqn: CRvsAB}.
\qed

\bigskip

Observe that  the normalization procedure transforms the solutions to the original problem  \eqref{primal} when this system is feasible.  More precisely, observe that $Ax = 0, \; x \ge 0$ if and only if $\tilde A D^{-1} x = 0, \; D^{-1} x \ge 0,$ where $D$ is the diagonal matrix whose diagonal entries are the norms of the columns of $A$.  Thus a solution to the original problem \eqref{primal}  can be readily obtained from a solution to the column-normalized preconditioned problem: premultiply by $D$.  At the same time, observe that the solution to \eqref{dual} does not change when the columns of $A$ are normalized.

Similarly, the balancing procedure transforms the solutions to the original problem \eqref{dual}.  In this case $A\transp y \ge 0, \; y\ne 0$ if and only if $\tilde A\transp R^{-1}y \ge 0, \; R^{-1}y\ne 0$, where $QR = A\transp$ is the QR-decomposition of $A\transp$.  Hence a solution to the original problem \eqref{dual} can be readily obtained from a solution to the row-balanced preconditioned problem: premultiply by $R$.  At the same time, observe that the solution to \eqref{primal} does not change when the rows of $A$ are balanced.

\bigskip

We next present our main result concerning properties of the two combinations of normalization and balancing procedures.

\begin{theorem}\label{main.thm} Assume $A = \matr{a_1 & \cdots & a_n}\in\R^{m\times n}$ is full-rank and $a_i\ne 0$ for $i=1,\dots,n$. 
\begin{description}
\item[(a)] Let $\hat A$ be the matrix obtained after normalizing the columns of $A$ and then balancing the rows of the resulting matrix.  Then
\begin{equation}\label{main.gcc}
\frac{1}{\sqrt{n}}\Cond_{GCC}(\hat A) \le  \Cond^{2,2}_R(\hat A) = \Cond_{Gr}(\hat A) \le\sqrt{n}\Cond_{GCC}(A).
\end{equation}
\item[(b)] Let $\hat A$ be the matrix obtained after balancing the rows of $A$ and then normalizing the columns of the resulting matrix.  Then
\begin{equation}\label{main.ab}
\frac{1}{\sqrt{n}} \Cond_{Gr}(\hat A) \le \Cond_R^{1,2}(\hat A) = \Cond_{GCC}(\hat A) \le\sqrt{n} \Cond_{Gr}(A).
\end{equation}
\end{description}
\end{theorem}
\proof
\begin{description}
\item[(a)] Let $\tilde A$ be the matrix obtained by normalizing the columns of $A$.   Proposition~\ref{prop.precond}(a) yields
\begin{equation}\label{eqn.1}
\Cond^{1,2}_R(\tilde A) =  \Cond_{GCC}(\tilde A) = \Cond_{GCC}(A).
\end{equation}
Since $\hat A$ is obtained by balancing the rows of $\tilde A$, Proposition~\ref{prop.precond}(b) yields
\begin{equation}\label{eqn.2}
\Cond^{2,2}_R(\hat A)= \Cond_{Gr}(\hat A) = \Cond_{Gr}(\tilde A) \le \Cond^{2,2}_R(\tilde A).
\end{equation}
Inequality \eqref{main.gcc} now follows by combining \eqref{ineq.norms}, \eqref{eqn: CRvsGCC}, \eqref{eqn.1}, and \eqref{eqn.2}.
\item[(b)] The proof this part is essentially identical to that of part (a) after using the parallel roles of the pairs $\Cond_{GCC}, \Cond^{1,2}_R$ and $\Cond_{Gr}, \Cond^{2,2}_R$  apparent from \eqref{eqn: CRvsGCC}, \eqref{eqn: CRvsAB}, and Proposition~\ref{prop.precond}.
\end{description}
\qed

Observe that both combined preconditioners in Theorem~\ref{main.thm} transform $A$ to a $\hat A = P A D$ where $D$ is a diagonal matrix and $P$ is a square non-singular matrix.  The particular $P$ and $D$ depend on what procedure is applied first.  Hence each of the combined preconditioners transforms both sets of solutions to the original pair of problems \eqref{primal} and \eqref{dual}.  In this case, solutions to the original problems can be readily obtained from solutions to the preconditioned problems via premultiplication by $D$ for the solutions to \eqref{primal} and by $P\transp$ for the solutions to \eqref{dual}.  

\medskip

As discussed in \cite[Section 5]{amelunxen2012coordinate}, there is no general relationship between $\Cond_{Gr}$ and $\Cond_{GCC}$, meaning that either one can be much larger than the other in a dimension-independent fashion.  Theorem~\ref{main.thm}(a) shows that when $\Cond_{GCC}(A) \ll \Cond_{Gr}(A)$  column normalization followed by row balancing would reduce both $\Cond_{Gr}$ and $\Cond^{1,2}_R$ while not increasing $\Cond_{GCC}$ beyond a factor of $\sqrt{n}$.  Likewise, Theorem~\ref{main.thm}(b) shows that when $\Cond_{Gr}(A) \ll \Cond_{GCC}(A)$  row balancing followed by column normalization would reduce both $\Cond_{GCC}$ and $\Cond^{1,2}_R$ while not increasing $\Cond_{Gr}$ beyond a factor of $\sqrt{n}$.  In other words, one of the two combinations of column normalization and by row balancing will always improve the values of all three condition numbers $\Cond_{GCC}(A)$, $\Cond_{Gr}(A)$ and $\Cond_R(A)$ modulo a factor of $\sqrt{n}$.
The following two questions concerning a potential strengthening of Thorem~\ref{main.thm} naturally arise:
\begin{itemize}
\item[(i)] Does the order of row balancing and column normalization matter in each part of Theorem~\ref{main.thm}? 
\item[(ii)] Are the leftmost and rightmost expressions in \eqref{main.gcc} and \eqref{main.ab} within a small power of $n$?  In other words, do the combined pre-conditioners make all $\Cond_{GCC}(A), \Cond_{Gr}(A)$ and $\Cond_R(A)$ the same modulo a small power of $n$?
\end{itemize}

The examples below, adapted from~\cite[Section 5]{amelunxen2012coordinate}, show that the answers to these questions are `yes' and `no' respectively.  Hence without further assumptions, the statement of Theorem~\ref{main.thm} cannot be  strengthened along these lines.

\begin{example}\label{exc order: bal&nor}
{\em 
Assume $\epsilon > 0$ and let $A = \matr{2/\epsilon & 1 & 1\\ 0 & -1  & 1}$. It is easy to see that $\Cond_{GCC}(A)=\sqrt{2}$.
After balancing the rows of $A$ and normalizing the columns of the resulting matrix we get \[ \hat A = \dfrac{1}{\sqrt{2(1+\epsilon^2)}}\matr{\sqrt{2(1+\epsilon^2)} & \epsilon & \epsilon \\[1.5ex] 0 &-\sqrt{2+\epsilon^2}&\sqrt{2+\epsilon^2} }.\]
It is easy to show that $\Cond_{GCC}(\hat A) = \frac{\sqrt{2\left(1+\epsilon^2\right)}}{\epsilon}$. Thus for $\epsilon > 0$ sufficiently small, $\Cond_{GCC}(\hat A)$  can be arbitrarily  larger than $\Cond_{GCC}(A)$.  Therefore  \eqref{main.gcc} does not hold for $A,\hat A$. }
\end{example}

\begin{example}{\em 
Assume  $0<\epsilon <1/2$ and let $A = \matr{-\epsilon & -1 & 1\\ 0 & -1  & 1+\epsilon}$.  Using \cite[Proposition 1.6]{amelunxen2012coordinate}, it can be shown that $\Cond_{Gr}(A) = \sqrt{2+\left(1+\epsilon\right)^2}$ . After normalizing the columns of $A$ and balancing  the rows of the resulting matrix, we obtain 
\[ 
\hat A = \frac{1}{\sqrt{\delta}}\matr{-\theta\delta & 
-\sqrt{2}\theta\epsilon\left(1+\epsilon\right)& -\theta\epsilon\sqrt{1+\left(1+\epsilon\right)^2} 
\\[1.5ex] 0 &-\sqrt{1+\left(1+\epsilon\right)^2}&\sqrt{2}\left(1+\epsilon\right)}. \]
 where $\delta = 1+3\left(1+\epsilon\right)^2$ and $\theta = \frac{1}{\sqrt{\delta+\epsilon^2}}$.  Using \cite[Proposition 1.6]{amelunxen2012coordinate} again, it can be shown that $\Cond_{Gr}(\hat A) = \sqrt{1+\frac{\delta}{\epsilon^2}} $.
Therefore, $\Cond_{Gr}(\hat A)$ can be arbitrarily larger than $\Cond_{Gr}(A)$. Hence \eqref{main.ab} does not hold for $A,\hat A.$ 
}
\end{example}

\begin{example}\label{largeGC} {\em 
Assume $\epsilon > 0$ and let $A = \matr{1+\epsilon & 1+\epsilon & -1+\epsilon\\ -1 & -1  & 1}$. In this case $\Cond_{GCC}(A) = \frac{\sqrt{2+2\left(1+\epsilon\right)^2} }{\epsilon}.$ 
After normalizing the columns of $A$ and balancing the rows of the resulting matrix, we get 
\[\hat A  = \dfrac{1}{\sqrt{2\gamma^2+2\beta^2}}\matr{ \beta& \beta& \sqrt{2}\gamma \\ -\gamma &-\gamma &\sqrt{2}\beta},\] 
where $\beta = \sqrt{\frac{1+\left(1+\epsilon\right)^2}{2}}$ and $\gamma= \sqrt{1+ \left(1-\epsilon\right)^2}$. It is easy to see that $\Cond_{GCC}(\hat A) = \sqrt{2}$. Therefore, $\Cond_{GCC}(\hat A)$ can be arbitrarily smaller than $\Cond_{GCC}(A)$ in \eqref{main.gcc}.}
\end{example}
\begin{example}{\em 
Assume $0<\epsilon < 1/2$ and let $A = \matr{2\epsilon & 1 & 1\\ 0 & -1  & 1}$. In this case $\Cond_{Gr}(A)  = \sqrt{1+\frac{1}{2\epsilon^2}}$. After balancing  the rows of  $A$ and normalizing the columns of the resulting matrix we get \[ \hat A = \dfrac{1}{\sqrt{2(1+\epsilon^2)}} \matr{\sqrt{2(1+\epsilon^2)} & 1& 1 \\[1.5ex] 0 &-\sqrt{1+2\epsilon^2}&\sqrt{1+2\epsilon^2} }.\]
Using \cite[Proposition 1.6]{amelunxen2012coordinate}, it can be shown that $\Cond_{Gr}(\hat A) = \sqrt{2+\epsilon^2}$. Thus  $\Cond_{Gr}(\hat A)$ can be arbitrarily smaller than $\Cond_{Gr}(A)$ in \eqref{main.ab}.
}
\end{example}

\bibliographystyle{plain}
\bibliography{mybibliography}

\end{document}